\documentclass[11pt]{article}

\usepackage{graphicx}
\usepackage{amsmath}
\usepackage{amssymb}

\newtheorem{thm}{Theorem}[section]
\newtheorem{dfn}{Definition}[section]
\newtheorem{cor}{Corollary}[section]
\newtheorem{lem}{Lemma}[section]
 	
\begin{document}

\title{
Alternating sign matrices with one $-1$ under vertical reflection}
\author{Pierre Lalonde 
\thanks{Supported by a grant from NSERC.}
\thanks{Mailing address: LaCIM, UQ\`AM, C.P. 8888, Succ.``A'',
 Montr\'eal, Qc, Canada, H3C 3P8; \newline
e-mail: lalonde@math.uqam.ca \,.} \\ 
LaCIM, Universit\'e du Qu\'ebec \`a Montr\'eal \\}
\date{}

\maketitle

\begin{abstract} 

We define a bijection that transforms an alternating sign matrix $A$ with 
one $-1$  into a pair $(N,E)$ where $N$ is a (so called) {\em neutral} alternating sign matrix (with one $-1$) and $E$ is an integer.  The bijection preserves the classical parameters of Mills, Robbins and Rumsey as well as three new parameters (including $E$).  It translates vertical reflection of $A$ into vertical reflection of $N$.  A hidden symmetry allows the interchange of $E$ with one of the remaining two new parameters.  
A second bijection transforms $(N,E)$ into a configuration of lattice paths called ``mixed configuration''.  
\end{abstract}

\vspace{2ex}
 
\section{Alternating sign matrices}  \label{class}

Recall that a square matrix $A=(a_{i j})_{1\leq i,j \leq n}$ is an 
order $n$ {\em alternating sign matrix} if $a_{i j}\in \{1,0,-1\}$ 
and if, in each row and each column, the non-zero entries alternate
in sign, beginning and ending with a $1$.  Thus, the entries of each row  and of each column add up to $1$.  

The entries in the first row of an alternating sign matrix are all $0$ 
except for one, which must be a $1$. It will be called the {\em 
first $1$}.

\vspace{2ex}

In their paper~\cite{MRR}, Mills, Robbins and Rumsey defined the 
following parameters on order $n$ alternating sign matrices $A=(a_{i j})$: 
\begin{itemize} 
  \item $r(A)$ is the number of entries to the left of the first $1$.  We 
         have $0\leq r(A) \leq n-1$.

  \item $s(A)$ is the number of entries that are equal to $-1$.

  \item $i(A)=\sum_{k>i,\ell<j}a_{i j}a_{k \ell}=
              \sum_{i,j}a_{i j}\left(\sum_{k>i,\ell<j}a_{k \ell}\right)$
              is the {\em number of inversions} of $A$. If $A$ is a permutation 
              matrix, $i(A)$ reduces to the usual number of inversions. 
\end{itemize}

We will use the following notation: $\mathcal{A}_{n}$ denotes the set of order $n$ alternating sign 
matrices and  $\mathcal{A}_{n,s}$ the set of order $n$ alternating sign 
matrices $A$ with $s(A)=s$. 

\vspace{2ex}

One of the Mills, Robbins and Rumsey conjectures asserts that $|\mathcal{A}_{n}|$ is also the number of order $n$ descending plane partitions.  In this form, the conjecture was solved by Zeilberger 
(see~\cite{Ze1},~\cite{Ze2}) with  subsequent simplifications by Kuperberg (see~\cite{Ku}).  Bressoud (see~\cite{Br})  gives an historical and mathematical account of the whole subject.

Stronger forms of the conjectures involve the parameters (defined above), which should translate into known combinatorially significant parameters on descending plane partitions.  In that direction, only special cases of the conjectures are solved. This is well known, of course, for $\mathcal{A}_{n,0}$ (permutation matrices).  The conjectures are also true for $\mathcal{A}_{n,1}$ (see~\cite{La1}).  This was done by encoding 
descending plane partitions into configurations of non-intersecting 
paths (so called TB-configurations),  which allows enumeration by a determinant.  After application of an algebraic transformation, the determinant is reinterpreted as enumerating another kind of lattice paths (mixed configurations), the set of which follows the same recurrences that describe  $\mathcal{A}_{n,1}$.  

In the present paper, we will give a bijective version of the last step, transforming $A\in \mathcal{A}_{n,1}$ into a pair $(N,E)$, where $N\in \mathcal{A}_{n,1}$ is ``neutral'' (to be defined in the next section) and $E$ is an integer.  A second bijection will transform the pair $(N,E)$ into a mixed configuration $\Omega$.  The bijections translate the already defined parameters (as well as three new ones) in a way that is coherent with the Mills, Robbins and Rumsey conjectures.

\vspace{2ex}

Let $A=(a_{i j})_{1\leq i,j \leq n}\in \mathcal{A}_{n}$. 
We write $\overline{A}=(a_{i, n+1-j})_{1\leq i,j \leq 
n}$ to denote the matrix obtained from $A$ by vertical reflection.  The classical parameters $r$, $i$ and $s$ applied to $A$ and to $\overline{A}$ are easily related (see~\cite{MRR}):
\begin{itemize}
    \item  $r(A)+r(\overline{A})=n-1$,
    \item  $i(A)+i(\overline{A})=\binom{n}{2}+s(A)$,
    \item  $s(\overline{A})=s(A)$.
\end{itemize}

Vertical reflection can be included in the conjectures.  It is then believed to correspond to an operation that can be interpreted as a kind of ``complementation'' operation on descending plane partitions.  In~\cite{La2}, it is shown that this operation takes a simple form in terms of  Gessel-Viennot paths duality (see~\cite{GV}) on TB-configurations.  (Krattenthaler (see~\cite{Kr}) has an even simpler interpretation in terms of rhombus tilings.) Our bijections  behave similarly: 
if $A\in \mathcal{A}_{n,1}$ is sent to $(N,E)$ and then to the mixed configuration $\Omega$, then $\overline{A}$ is sent to $(\overline{N},-E)$, which is sent to $\overline{\Omega}$, the Gessel-Viennot dual of $\Omega$.

\vspace{3ex}

\section{Three new parameters} \label{new}

\begin{figure}
    \centerline {
      \includegraphics[scale=1.0]{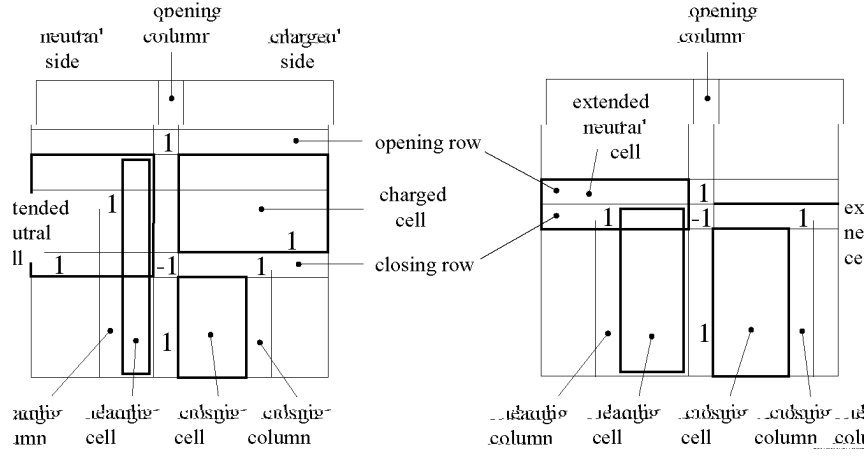}
     }
    \caption{Schematic view of a positive matrix (left) and a neutral one (right), with some of the related regions as defined in this section.  Only the significant non-zero entries are depicted.  (The {\bf {\large 0}} region contains only 0's.) \label{fig1}}
\end{figure}

In what follows, we will introduce the three new parameters defined for a matrix $A\in \mathcal{A}_{n,1}$.  These parameters are related to various sub-matrices of $A$, which we describe below (see also figure~\ref{fig1}). 

\begin{itemize} 
  \item The {\em opening column} of $A$ is the column of its (unique)  $-1$.  The highest $1$ in this column is the {\em opening $1$} and 
  the corresponding row, the {\em opening row}.  The {\em closing row} is the row of the $-1$. The opening column divides $A$ into a {\em left side} and a {\em right side} (both excluding the opening column).
  
  \item   The closing row is the only row that contains two $1$, one in each side.   These $1$ will be referred to as the {\em left} $1$ and the {\em right} $1$.
  
  \item If any, the rows between the opening and the closing rows 
  are the {\em enclosed rows}.  If there are no enclosed rows, $A$ is said to be {\em neutral}; otherwise $A$ is {\em charged}. In the latter case, define the {\em charged side} to be the side (left or right) where we find the  $1$ of the lowest enclosed row, the other side being the {\em neutral side}.  If the charged side is the right side (respectively: left side), we say that $A$ is {\em positive} (respectively: {\em negative}).  
    
\end{itemize}

(In fact, we can define more generally $A=(a_{i j})\in \mathcal{A}_{n}$ to be {\em neutral} if $a_{i j}=1$ when $a_{i+1\,j}=-1$.) 

\vspace{2ex}

Let $\mathcal{A}_{n,1}^+$ (respectively: $\mathcal{A}_{n,1}^0$, $\mathcal{A}_{n,1}^-$) be the set of positive  (respectively: neutral, negative) matrices $A\in \mathcal{A}_{n,1}$.  These sets are mutually disjoint and form a partition of $\mathcal{A}_{n,1}$.  Moreover, $\mathcal{A}_{n,1}^{+}$ and $\mathcal{A}_{n,1}^{-}$ are mirror-images of one another: $A\in \mathcal{A}_{n,1}^{+}$ iff $\overline{A}\in \mathcal{A}_{n,1}^{-}$.

  We further define the following for $A\in \mathcal{A}_{n,1}^+\cup \mathcal{A}_{n,1}^0$:
\begin{itemize}
  \item The intersection of the enclosed rows with the right (respectively: left) side  defines the {\em charged} (respectively: {\em neutral}) {\em cell}.  The {\em extended} neutral cell includes the intersection of the opening and of the closing rows with the left side.  If $A\in \mathcal{A}_{n,1}^0$,  the charged and the neutral cells are empty.
  
  \item The highest $1$ in the left side below the opening row is the {\em leading $1$}.  Its column is the {\em leading column}.  The sub-matrix between the leading and the opening column and below the opening row is the {\em leading cell}.  The sum of the entries of the leading cell is denoted $\ell(A)$.
    
  \item Finally,  the right $1$ (in the closing row) is also called the {\em closing $1$}.  Its column is the {\em closing column}. The sub-matrix of $A$ between  the closing and the opening column and below the closing row is the {\em closing cell}. The {\em extended} closing cell includes the parts of opening and of the closing columns that are below the closing row.  The sum of the entries of the closing cell is denoted $c(A)$.

\end{itemize}
  
{\bf Remark.\,\,} It should be observed that $\ell(\overline{A})=c(A)$ and  $c(\overline{A})=\ell(A)$ when $A\in \mathcal{A}_{n,1}^{0}$.

\vspace{2ex}

\begin{figure} 
    \centerline {
      \includegraphics[scale=1]{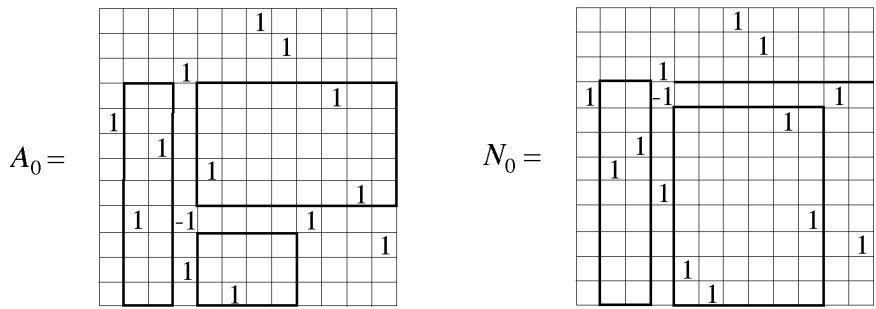}
     }
    \caption{In each of the above matrices, the leading, charged and closing cells are emphasized.   Matrix  $A_{0}$ is positive, with $E(A_{0})=3$, $B(A_{0})=-1$ and $J(A_{0})=7$. Matrix $N_{0}$ is neutral, with $E(N_{0})=0$, $B(N_{0})=2$ and $J(N_{0})=7$. The classical parameters are: $r(A_{0})=r(N_{0})=6$ and $i(A_{0})=i(N_{0})= 30$. \label{fig2}}
\end{figure}

We can now define the new parameters (see figure~\ref{fig2}):
\begin{itemize} 
  
  \item  If $A\in \mathcal{A}_{n,1}^{+}$, its {\em electric charge}, $E(A)$, is the sum of the entries of the charged cell of $A$. In that case, $E(A)>0$.  Define $E(A)=0$ if $A\in \mathcal{A}_{n,1}^{0}$ and $E(A)=-E(\overline{A})$ if $A\in \mathcal{A}_{n,1}^{-}$.  Thus $A$ is positive, neutral or negative according to the sign of $E(A)$.

  \item If $A\in \mathcal{A}_{n,1}^{+}\cup  \mathcal{A}_{n,1}^{0}$, define its {\em magnetic charge} by $B(A)=c(A)-\ell(A)$.   If $A\in \mathcal{A}_{n,1}^{0}$, we clearly have $B(\overline{A})=-B(A)$.  Extend this property to define $B(A)$ for $A\in  \mathcal{A}_{n,1}^{-}$.
  
    \item If $A\in \mathcal{A}_{n,1}^{+}\cup  \mathcal{A}_{n,1}^{0}$, define $J(A)=c(A)+\ell(A)+|E(A)|+1$.  Notice that $J(A)=J(\overline{A})$ if  $A\in \mathcal{A}_{n,1}^{0}$.   Extend this property to define $J(A)$ for $A\in  \mathcal{A}_{n,1}^{-}$.
\end{itemize}
\vspace{2ex}

\vspace{2ex}

Clearly, with respect to vertical reflection, $E$ and $B$ are anti-invariants, while $J$  is  invariant. Algebraically: 
$$E(A)+E(\overline{A})=0, \qquad B(A)+B(\overline{A})=0 \qquad \mbox{and} \qquad J(\overline{A})=J(A)\,.$$ 

\vspace{3ex}

\section{Neutralizing alternating sign matrices}

Our first task will be to learn how to ``neutralize'' a given matrix $A\in \mathcal{A}_{n,1}^{+}$.  This requires many steps based on the {\em horizontal/vertical displacement} procedure.  

\vspace{2ex}

{\bf Horizontal displacement} ($H$):
Let $P=(p_{i j})_{1\leq i \leq m, 1\leq j \leq n}$ be a $(0,1)$-matrix.  Suppose that the non-zero columns occupy positions $j_{1}<j_{2}< \cdots < j_{k}$ with  $j_{1}=1$ and $j_{k}<n$.

Its horizontal displacement, $H(P)$, is the matrix obtained from $P$ by displacing the entries of column $j_{i}$ to column $j_{i+1}$ (for $1\leq i \leq k$), where $j_{k+1}=n$.  Column $j_{1}=1$ is replaced by a column of $0$'s.  Clearly, $H(P)$ is a $(0,1)$-matrix of the same dimension as $P$, with non-zero columns in positions  $j_{2}<\cdots j_{k}<j_{k+1}=n$.  The procedure is obviously injective.

\vspace{2ex}

We define similarly the {\em vertical displacement} $V(P)$ for $(0,1)$-matrices $P$ such that  the first row is $0$ and the last, non-zero. (The rows are displaced from bottom to top.) 

\vspace{2ex} 

   For instance, 
$$V\begin{pmatrix}
     0 & 0 & 0 & 0  \\
     0 & 0 & 0 & 1  \\
     0 & 0 & 0 & 0  \\
     0 & 0 & 0 & 0  \\
     0 & 0 & 1 & 0  \\
     1 & 0 & 0 & 0    
      \end{pmatrix}=
     \begin{pmatrix}
     0 & 0 & 0 & 1  \\
     0 & 0 & 1 & 0  \\
     0 & 0 & 0 & 0  \\
     0 & 0 & 0 & 0  \\
     1 & 0 & 0 & 0  \\  
     0 & 0 & 0 & 0  
      \end{pmatrix}\,.$$

\vspace{2ex}

We will apply the horizontal/vertical displacement to some of the cells of a given matrix $A\in \mathcal{A}_{n,1}^{+}\cup \mathcal{A}_{n,1}^{0}$ (or to some modifications of $A$).  This will give the {\em discharging procedure} which essentially transforms $A$ into a permutation matrix $P$ of the same dimension.  In what follows, the opening column, closing cell,\ldots of any transformation of $A$ refer to sub-matrices of the transformed matrix that occupies the same position as in $A$.

\vspace{2ex}

\begin{dfn}
{\bf (Partial discharging procedure $\delta$)}
 Let $A\in  \mathcal{A}_{n,1}^{+}\cup \mathcal{A}_{n,1}^{0}$. The corresponding discharged matrix $\delta(A)$ is obtained by successively applying the following rules to $A$: 
 \begin{enumerate}
  \item Erase the $-1$ and the closing $1$.  \label{s1}
  \item Apply $H$ to the extended closing cell.  \label{s2}
  \item Apply $V$ to the extended neutral cell.  \label{s3}
  \item Lower the $1$'s in the extended neutral and in the charged cells by one row (erasing or writing $0$'s when necessary). \label{s4}
  \end{enumerate}
\end{dfn}


If $A$ is neutral, observe that step 4 cancels the effect of step 3.  Hence we only  need to apply steps 1 and 2.

\begin{figure} 
    \centerline {
      \includegraphics[scale=1.0]{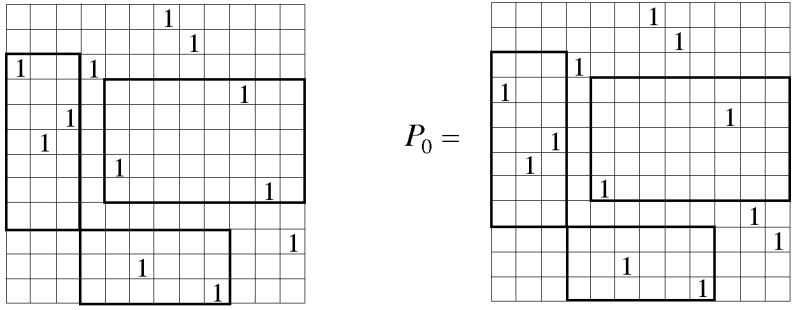}
     }
    \caption{ The first matrix results from $A_{0}$ (figure~\ref{fig2}) after applying the first three steps of $\delta$.  The last matrix is $P_{0}=\delta(A_{0})$.  The emphasized sub-matrices are the extended neutral cells, the charged cells and the extended closing cells.   \label{fig3}}
\end{figure}

\begin{lem} \label{bijmoins}
Let $A\in \mathcal{A}_{n,1}^{+}\cup \mathcal{A}_{n,1}^{0}$. Let $k$ be the position of its opening row and $P=(p_{i j})=\delta(A)$.  Then
\begin{enumerate}
  \item $P$ is a permutation matrix.
  \item rows 1 to $k$ (included) are the same in $A$ as in $P$.  Thus $r(A)=r(P)$.
  \item if $p_{k j}=1$ and $p_{k+1, m}=1$ then $m<j$.
  \item $\ell(A)=\ell(P)$, where $\ell(P)$ is the sum of the entries of $P$ in its leading cell (i.e. below row $k$ and strictly between columns $m$ and $j$).  
  \item $c(A)+E(A)< x(A)=x(P)$.  (Here,  $x(A)$ (respectively: $x(P)$) denotes the sum of the elements of $A$ (respectively: of $P$) that are below the opening row (row $k$) and in the right side.)
  \item $i(A)=i(P)+c(A)+1+E(A)$.
\end{enumerate}
\end{lem}

\vspace{2ex}

{\bf Proof.\,\,}  Let $A\in \mathcal{A}_{n,1}^{+}\cup \mathcal{A}_{n,1}^{0}$.   Notice that, after step~\ref{s1}, the resulting matrix is a $(0,1)$-matrix that differs from a permutation matrix only in the opening column (two $1$'s, one being below the closing row) and in the closing column (no $1$'s).  Thus we can apply step~\ref{s2}, resulting in a permutation matrix.  The left part of the closing row contains the left 1, allowing the application of step~\ref{s3}.  Now the opening row contains two $1$ (one in the left side), the closing row none.  The other rows and the columns contain one $1$.  (See for instance,  the first matrix of figure~\ref{fig3}).   Clearly, we can then apply step 4, to get a permutation matrix  (e.g., the last matrix of figure~\ref{fig3}).

The rows from the first to the opening row are unaffected by $\delta$.  The leading $1$ of $A$ is vertically displaced to the row just below the opening row (proving the third statement).  All of the $1$'s in the leading cell of $A$ remain in the leading cell of $P$, showing that $\ell(P)=\ell(A)$. 

   Applying the discharging procedure, the region accounted for by $x(A)$ loses a $1$ (the closing $1$), but gains one $1$ (the lowest $1$ of the opening column, after step 2).  Thus $x(P)=x(A)>c(A)+E(A)$.
 
Let us examine the behavior of the number of inversions.  Let $y$ be the sum of the entries of $A$ that are strictly between the opening and the closing columns and are strictly higher than the closing row.  Apply $\delta$ to $A$.  Step 1, the erasure of the $-1$ and of the closing $1$, deletes $c(A)+1-y$ inversions (related to the elements in the intersection of the South-West and of the North-East regions of each of the erased elements).  Step 2 deletes $y$ inversions (to see this, observe that any two $1$'s, not both between the opening and the closing columns, keep their relative positions).  Similarly, steps 3 and 4 (combined) delete $E(A)$ inversions, giving a total of $c(A)+1+E(A)$ deleted inversions.  Thus $i(A)=i(P)+c(A)+1+E(A)$. 
\quad $\Box$

\vspace{2ex}

Next, we complete $\delta$ with the necessary information to get a bijection.
    
\vspace{2ex}

\begin{dfn}
{\bf (Complete discharging procedure $\Delta$)}
Let $A\in \mathcal{A}_{n,1}^{+}\cup \mathcal{A}_{n,1}^{0}$.  We define
$$\Delta(A)  =  (k,\delta(A), c(A),E(A))\,,$$ where $k$ is the position of the opening row of $A$ ($1\leq k\leq n$).  
\end{dfn} 

\vspace{2ex}

{\bf Examples.\,\,} Referring to figures~\ref{fig2} and~\ref{fig3}, we have $\Delta(A_{0})=(3,P_{0},1,3)$ and $\Delta(N_{0})=(3,P_{0},4,0)$.

\vspace{2ex}

We now proceed to determine the range of $\Delta$.
Let $\mathcal{B}_{n,1}$ be the set of 4-tuples $(k,P,c,E)$ such that:
\begin{enumerate}
  \item $1\leq k \leq n-2$,
  \item $P=(p_{i j})\in \mathcal{A}_{n,0}$,
  \item if $p_{k j}=1$ and $p_{k+1, m}=1$ then $m<j$,
  \item $c\geq 0$, $E\geq 0$ and $c+E< x(P)$. 
\end{enumerate}

\vspace{2ex}

\begin{lem} \label{bijplus}
The discharging procedure $\Delta$ is a bijection from  $\mathcal{A}_{n,1}^{+}\cup \mathcal{A}_{n,1}^{0}$ to $\mathcal{B}_{n,1}$.
\end{lem}

{\bf Proof.\,\,}  By lemma~\ref{bijmoins}, it is clear that $\Delta(\mathcal{A}^{+}_{n,1}\cup \mathcal{A}_{n,1}^{0}) \subseteq \mathcal{B}_{n,1}$.  Let $(k,P,c,E) \in \mathcal{B}_{n,1}$.  Suppose that we can find the location of the charged, extended neutral and extended closing cells in $P$ from the given information, then we can readily reverse each step of $\delta$, showing $\Delta$ to be invertible.   This is easy: $k$ locates in $P$ the opening row.  Since $P$ is a permutation matrix, there is a unique $1$ in this row, which defines the opening column, which in turn defines the right and the left sides. 

The closing row is the highest row below the opening row such that the elements between (and including) these rows  in the right side sum up to $E$.  Thus we can find the charged cell and the extended neutral cell.  Since the (unique) $1$ in row $k+1$ is in the left side, we can apply the reverse of step 4 and  the reverse of step 3 of $\delta$.  If $E>0$, the closing row had a $1$ in the right side, which now  ends in the charged cell (hence the resulting matrix will be positive), emptying the right part of the closing row. If $E=0$, the right part of the closing row (position $k+1$) was already empty (and the reverses of steps 4 and 3 cancel each other).
 
 The closing column is the leftmost column to the right of the opening column such that the elements between (and including) these columns and below (strictly) the closing row sum up to $c+1$.  This defines the extended closing cell.  The construction is always possible since $c+1+E\leq x(P)$. Notice that the rightmost column of this cell contains a $1$ and its leftmost is empty (allowing to apply the reverse of step 2).  It remains to place the $-1$ and the closing $1$; this is readily done since the right part (including the opening column) of the closing row contains no $1$.
 \quad $\Box$

\vspace{2ex}

Now comes the goal of this section, the definition of the neutralizing procedure.  It is based on the following remark:  let $c,E\geq 0$ then: $(k,P,c,E) \in \mathcal{B}_{n,1}$ iff $(k,P,c+E,0) \in \mathcal{B}_{n,1}$.  Moreover $\Delta^{-1}(k, P, c+E,0)\in \mathcal{A}_{n,1}^{0}$.

\vspace{2ex}

\begin{dfn}{\bf (Neutralizing procedure $\Lambda$, non-negative case)} 

Let $A\in  \mathcal{A}_{n,1}^{+}\cup \mathcal{A}_{n,1}^{0}$ with $(k,P,c,E)=\Delta(A)$. We define 
$\Lambda(A)=(\Delta^{-1}(k, P, c+E,0), E)$.
\end{dfn}

Observe that if $A$ is neutral then $E(A)=0$, leading to $\Lambda(A)=(\Delta^{-1}(k, P, c+0,0), 0)=(A,0)$.  Since $\overline{A}$ also is  neutral, we have $\Lambda(\overline{A})=(\overline{A},0)$, a property that we will extend to define $\Lambda$ over $\mathcal{A}^{-}_{n,1}$.

\begin{dfn}{\bf (Neutralizing procedure $\Lambda$, negative case)} 

\vspace{1ex}
Let $A\in  \mathcal{A}_{n,1}^{-}$ (so that $\overline{A}\in  \mathcal{A}_{n,1}^{+}$).  Writing $(\overline{N},-E)= \Lambda(\overline{A})$, we define $\Lambda(A)=(N,E)$.
\end{dfn}

\vspace{2ex}

\begin{figure} 
    \centerline {
     \includegraphics[scale=1]{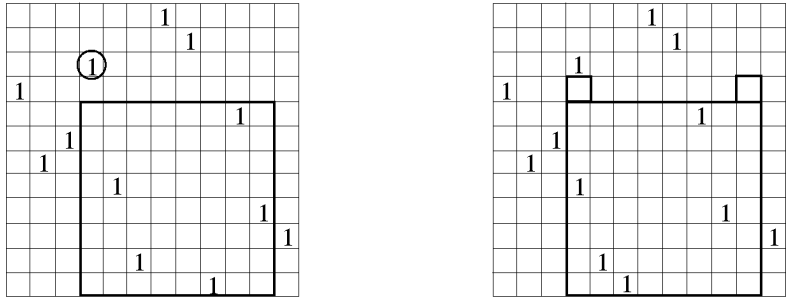}
     }
    \caption{Determination of $N_{0}=\Delta^{-1}(3, P_{0}, 4,0)$ for matrix $P_{0}$, the last matrix of figure~\ref{fig3}. \label{fig4}}
\end{figure}

{\bf Example.\,\,}  We will construct $\Lambda(A_{0})$ (from figures~\ref{fig2} and~\ref{fig3}).   We already know that $\Delta(A_{0})=(3,P_{0},1,3)$ so that $\Lambda(A_{0})=(N_{0},3)$ where $N_{0}=\Delta^{-1}(3, P_{0}, 4,0)$ is the neutral matrix from figure~\ref{fig2}.   To see how to compute it from $P_{0}$, look at figure~\ref{fig4}: 
\begin{enumerate}
  \item The first matrix is  $P_{0}$.  Row 3 is the opening row and row 4 is the (new) closing row.  The position of the $1$ (circled) in row 3 defines the opening column (column 4).  Since $N_{0}$ is neutral, it suffice to find the extended closing cell (emphasized; it must contain five $1$'s).
  \item We then apply  $H^{-1}$ to the extended closing cell (reverse of step 2), giving the second matrix.  It remains to write the $-1$ in position $(4,4)$ and  the closing $1$ in position $(4,11)$ (small squares).
\end{enumerate}


\begin{thm} \label{clef}
The neutralizing procedure is a bijection
$$\Lambda: \mathcal{A}_{n,1} \longrightarrow \mathcal{N}_{n,1}:=\{(N,E)\,|\, N\in \mathcal{A}_{n,1}^0, E\in \mathbb{Z}, -\ell(N)\leq E \leq c(N)\}\,.$$

Moreover, let $A\in  \mathcal{A}_{n,1}$ and $\Lambda(A)=(N,E)$, then:
\begin{enumerate}
  \item $A\in \mathcal{A}_{n,1}^0$ iff $N=A$. \label{neutre}
  \item $\Lambda(\overline{A})=\overline{\Lambda(A)}:=(\overline{N},-E)$.\label{dualdelta}
  \item The position $k$ of the opening row is the same in $A$ as in $N$. 
  In fact, rows 1 to $k$ are the same in $A$ and in $N$.  Thus $r(A)=r(N)$.  \label{ligne}
  \item The following relations hold:   \label{param}  
  \begin{enumerate}
  \item $i(N)=i(A)$,
  \item $E=E(A)$,
  \item $B(N)=B(A)+E(A)$,
  \item $J(N)=J(A)$.
\end{enumerate}   
\end{enumerate}
\end{thm}

{\bf Proof.\,\,} Partition $\mathcal{N}_{n,1}$ into three subsets $\mathcal{N}_{n,1}^{+}$, $\mathcal{N}_{n,1}^{0}$ and $\mathcal{N}_{n,1}^{-}$, defined according to the sign ($+$, $0$ or $-$) of the second component ($E$).  Statement~\ref{neutre} (which is already known in one direction and trivial in the other) shows that $\Lambda$  
bijectively maps $\mathcal{A}_{n,1}^0$ to  $\mathcal{N}_{n,1}^{0}$.  The other statements trivially follow in the neutral case.

Suppose $A\in  \mathcal{A}_{n,1}^{+}$.  Let $(k,P,c,E)=\Delta(A)$ and $(N,E)=\Lambda(A)$.  Since $\Delta(N)=(k,P,c+E,0)$, we have $c(N)=c+E>0$, showing that $c(N)\geq E>0$   (so that $\Lambda(\mathcal{A}_{n,1}^{+})\subseteq \mathcal{N}_{n,1}^{+}$).  Conversely, let $(N,E)\in \mathcal{N}_{n,1}^{+}$, with $(k,P,c,0)=\Delta(N)$  (hence, $c\geq E>0$ and $c<x(P)$).  Thus $(k,P,c-E,E)\in \mathcal{B}_{n,1}$ and there is a unique $A\in \mathcal{A}_{n,1}^{+}$ such that $\Delta(A)=(k,P,c-E,E)$.  Consequently, $\Lambda$  
bijectively maps $\mathcal{A}_{n,1}^+$ into  $\mathcal{N}_{n,1}^{+}$.
Statement~\ref{ligne} clearly holds (since $\delta(A)=P=\delta(N)$). Statement~\ref{param}  follows easily from lemma~\ref{bijmoins} and the obvious relations $c(N)=c(A)+E(A)$ and $\ell(N)=\ell(A)$.

Statement~\ref{dualdelta} is already known when $A$ is neutral and is true by construction when $A$ is charged.   Using $\Lambda(A)=\overline{\Lambda(\overline{A})}$, we see that $\Lambda$ bijectively maps $\mathcal{A}_{n,1}^-$ to  $\mathcal{N}_{n,1}^{-}$.  Thus, if $A$ is negative, statement~\ref{ligne} holds by symmetry.  Statement~\ref{param} follows easily from the positive case, with the help of the formulae at the end of sections~\ref{class} and~\ref{new}.  \quad $\Box$

\vspace{3ex}

\section{Exchanging the electric charge and the magnetic charge}

Using the neutralizing procedure, we define an involution on $\mathcal{A}_{n,1}$ that exchanges $E$ and $B$.  Thus the two charges play the same r\^ole and are completely interchangeable.  We begin by showing  that both charges share a common range.

\begin{lem} \label{range}
Let $A\in \mathcal{A}_{n,1}$ and $(N,E)=\Lambda(A)$. Then $-\ell(N)\leq B(A) \leq c(N)$.
\end{lem}  

{\bf Proof.\,\,}  From theorem~\ref{clef}, we know that $B(A)=c(N)-\ell(N)-E(A)$.  Thus $B(A)=c(N)-(\ell(N)+E(A))\leq c(N)$, since $E(A)\geq -\ell(N)$.  On the other hand, we can also write $B(A)=(c(N)-E(A))-\ell(N)\geq -\ell(N)$, since $E(A)\leq c(N)$. \quad $\Box$ 

\vspace{2ex}

Observe that the function defined on $\mathcal{N}_{n,1}$ by $\xi(N,E)=(N,c(N)-\ell(N)-E)$ is in fact an involution.  Thus $\Lambda^{-1}\circ \xi \circ \Lambda$ is an involution on $\mathcal{A}_{n,1}$.  We will write $A'$ instead of the more cumbersome $(\Lambda^{-1}\circ \xi \circ \Lambda)(A)$.

\vspace{2ex}

\begin{thm}  Let  $A\in \mathcal{A}_{n,1}$.  We have:
\begin{enumerate}
  \item  The position $k$ of the opening row is the same in $A$ as in $A'$. 
  In fact, rows 1 to $k$ are the same in $A$ as in $A'$.  Thus $r(A)=r(A')$. \label{ex3}
  \item Parameters $i$ and $J$ take the same values on $A$ as on $A'$. \label{ex5}
  \item The involution $(\quad)'$ exchanges the charges: $E(A')=B(A)$ and $B(A')=E(A)$. \label{ex4}
  \item $(\,\overline{A}\,)'=\overline{(A')}$. \label{ex2}

\end{enumerate}
\end{thm}

{\bf Proof.\,\,}   Let $\Lambda(A)=(N,E)$ so that $\Lambda(A')=(N,c(N)-\ell(N)-E)$.  By theorem~\ref{clef},  statements~\ref{ex3} and~\ref{ex5} are true.  For the charges, we have $E(A')= c(N)-\ell(N)-E=B(A)$ and $B(A')=c(N)-\ell(N)-E(A')=E(A)$.

Finally, statement~\ref{ex2} is easily proved:
\begin{eqnarray*}
(\Lambda^{-1}\circ \xi \circ \Lambda)(\overline{A}) & = & 
                        (\Lambda^{-1}\circ \xi)(\overline{N},-E) \\
          & = & \Lambda^{-1}(\overline{N},c(\overline{N})-\ell(\overline{N})+E)\\
          & = & \Lambda^{-1}(\overline{N},-(c(N)-E-\ell(N))\\
          & = & \overline{A'}\,. \qquad \Box
\end{eqnarray*}
          
\vspace{2ex}

\begin{figure} 
    \centerline {
     \includegraphics[scale=1.0]{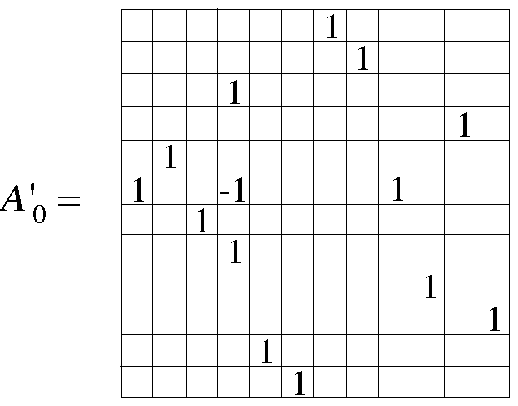}
     }
    \caption{
       Matrix $A'_{0}$, where $A_{0}$ is defined in figure~\ref{fig2}.   \label{fig5} }
\end{figure}

\vspace{2ex}

Of course, this leads to another bijection, $\xi\circ\Lambda: \mathcal{A}_{n,1} \longrightarrow \mathcal{N}_{n,1}$, which focuses on the parameter $B$ instead of $E$.  In fact, $(\xi\circ\Lambda)(A)=(N,B(A))$.

\vspace{3ex}

\section{Encoding elements of $\mathcal{N}_{n,1}$ into mixed configurations}

It is well known that a permutation matrix $P=(p_{i j})\in  \mathcal{A}_{n,0}$ can be bijectively encoded by a sequence $(a_{i})_{i=1}^{n}$ of non-negative integers called its inversion table.  In fact, $a_{i}$ is the sum of the entries of $P$ that are below row $n+1-i$ and to the left of the unique $1$ in that row.  With this convention, we have $0\leq a_{i}<i$ for $1\leq i\leq n$.  The classical parameters are easily recovered: $r(P)=a_{n}$ and $i(P)=a_{1}+\cdots+a_{n}$. Clearly, $(\overline{a}_{i})_{i=1}^{n}$, where  $\overline{a}_{i}=i-1-a_{i}$ (for $1\leq i\leq n$), is the inversion table of $\overline{P}$.  We will define a generalization of inversion tables that applies to $\mathcal{A}_{n,1}$.

\vspace{2ex}

\begin{dfn}
Let  $(N,E)\in \mathcal{N}_{n,1}$.  Let $n+1-k$ be the position of the opening row of $N$ (thus the position of the closing row is $n+2-k$).  For $1\leq i\leq n$, define $a_{i}$ as the sum of the entries of $N$ that are below row $n+1-i$ and to the left of the unique $1$ (or the leftmost $1$ if $i=k-1$) in that row.   Let $b=c(N)$ and $\beta=E+\ell(N)$.  The sequence of non-negative integers $(k;a_{1},\ldots,a_{n};b,\beta)$ is called the {\em generalized inversion table} of $(N,E)$.  
\end{dfn}

\vspace{2ex}

{\bf Remarks.\,\,} 
\begin{itemize}
\item The element $a_{i}$ is non-negative even if the sum (defining $a_{i}$) includes the $-1$, since it will include (at least) the left $1$.
\item By definition, $b=c(N)\geq 0$ and $\beta=E+\ell(N)\geq 0$. 
\item Clearly, $\ell(N)=a_{k}-1-a_{k-1}$, an observation that we will often use later.
\end{itemize}
 
\vspace{2ex}

{\bf Example.\,\,}  For instance, the generalized inversion table of $(N_{0},3)$ (from figure~\ref{fig2}) is: $$(10;\,\, 0,0,2,2,0,0,1,5,0,3,6,6;\,\, 4,5)\,.$$

\vspace{2ex}

 \begin{lem} \label{char}
A sequence $(k; a_{1},\ldots,a_{n};b,\beta)$
of non-negative integers is the generalized inversion table of  some unique $(N,E)\in \mathcal{N}_{n,1}$ iff
\begin{enumerate}
    \item $3\leq k \leq n$,
    \item  $0\leq a_{i}\leq i-1$ for all $i$,
    \item  $a_{k-1}<a_{k}$,
    \item  $a_{k-1}+\beta< a_{k}+b\leq k-2$.
\end{enumerate}
\end{lem}

{\bf Proof.\,\,}  First, we show that the conditions are necessary.  Clearly $1\leq n+1-k\leq n-2$, which is equivalent with  $3\leq k \leq n$.
 The elements below row $n+1-i$ of $N$ add up to $i-1$, forcing $a_{i}\leq i-1$.  Similarly, focusing on the elements below row $n+2-k$ leads to $a_{k}+b\leq k-2$   (in the sum defining $a_{k}$, replace the left $1$ by  the $1$ below the $-1$).   The fact that the leftmost $1$ in the closing row is in the left side translates into $a_{k-1}<a_{k}$.  Using $\beta=E+a_{k}-1-a_{k-1}$, we get $a_{k-1}+\beta=a_{k}+E-1<a_{k}+c(N)=a_{k}+b$. 

The conditions are sufficient: any sequence $(k; a_{1},\ldots,a_{n};b,\beta)$ verifying the conditions of the lemma is the inversion table of a unique element $(N,E)\in \mathcal{N}_{n,1}$.  In fact, as with ordinary inversion table, we can construct $N$ from the first row down.  Since $3\leq k \leq n$, the integer $k$ determine the opening row (and hence the closing row) of $N$.  The non-negative numbers $a_{i}$ (from $i=n$ down to $i=n+2-k$) are used to position the (unique, or leftmost) $1$'s in rows $1$ to $n+2-k$ (closing row) of $N$.  Now that we know the opening column, we can place the $-1$.   Since $a_{k-1}<a_{k}$, the leftmost $1$ in the closing row is in the left side of $N$.  To place the rightmost $1$ in the closing  row, use $b=c(N)$. (Since $a_{k}\leq a_{k}+b<k-2$, it is always possible to find the position of this $1$.) Place the remaining $1$'s in row $n+3-k$ to $n$ according to the remaining $a_{i}$'s.   This completely determine $N$.  Finally, define $E=a_{k-1}+1-a_{k}+\beta$.  Since $a_{k-1}+\beta< a_{k}+b$, we have $E\leq b=c(N)$.  Since $\ell(N)=a_{k}-1-a_{k-1}$, we also have $E=\beta-\ell(N)\geq -\ell(N)$.  \quad $\Box$
 
\vspace{2ex}

Inversion tables are encoded as sequences of non-intersecting lattice paths called {\em mixed configurations}.  We consider lattice-paths on the strict half-grid $\mathcal{G}_{n}=\{(k,\ell)\,|\,{0\leq  k < \ell \leq n}\}$.   (For more symmetry, the grid will be slightly shifted so that its boundary forms a reversed equilateral triangle.) {\em Mixed paths} on $\mathcal{G}_{n}$ are composed  of two consecutive parts (Left and Right), where:
\begin{itemize}
   \item the Left part is composed of South  steps (S) and East steps (E).  
   \item the Right part is composed of (another kind of) East 
   steps (F) and North-East steps (N). 
\end{itemize}
Notice that each path contains a vertex that belongs both to the Left part and to the Right part.  Such vertices are called {\em junctions}. 

\vspace{2ex}

An  order $n$ {\em mixed configuration} is a sequence of mixed paths
$\Omega =(\omega_1,\ldots,\omega_n)$ on $\mathcal{G}_{n}$ such that:
\begin{itemize}
	\item There is a permutation $\sigma$ such that $\omega_{i}$ starts     from $(0,i)$ and ends at $(\sigma(i)-1,\sigma(i))$. 
	\item The sub-configuration obtained by deleting the Right part (respectively: Left part) of paths is non-intersecting (no common vertex).
\end{itemize}

 We write $\mathcal{M}_{n,s}$ to denote the set of order $n$ mixed configurations with $s$ N-steps.

\vspace{2ex}

\begin{figure} 
	\begin{center}
           \includegraphics[scale=1.0]{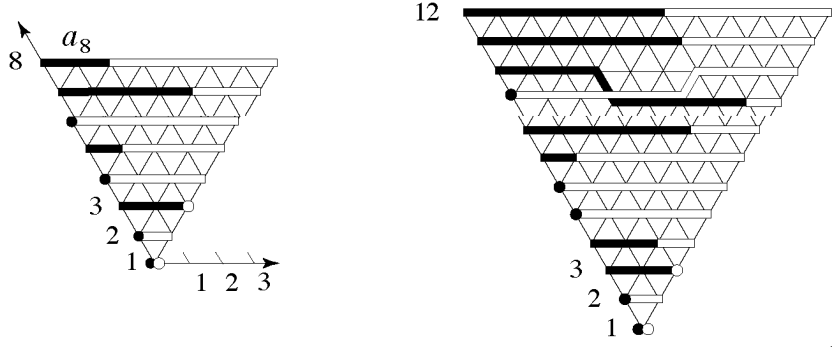} 
	\end{center}
	\caption{A mixed configuration with no N-step (left). A mixed configuration with one N-step (right). The Left part of each is colored black; the Right part, white. \label{fig6}}
	\label{perm}
\end{figure}

Observe that if $\Omega =(\omega_1,\ldots,\omega_n)\in \mathcal{M}_{n,0}$, then all paths of $\Omega$ are horizontal  (see figure~\ref{fig6} (left)).   The corresponding  inversion table is the sequence of the lengths of the Left parts of the paths.
 
\vspace{2ex}

More interesting are  configurations 
$\Omega =(\omega_1,\ldots,\omega_n)\in \mathcal{M}_{n,1}$. 
In that case, $\Omega$ has  two consecutive special paths 
$\omega_{k-1}$ (which contains the N-step) and $\omega_{k}$  (which contains a S-step). The other paths are horizontal (see figure~\ref{fig6} (right)).

\begin{dfn} \label{defphi}
Let $(N,E)\in \mathcal{N}_{n,1}$ and $(k;a_{1},\ldots,a_{n};b,\beta)$ its generalized inversion table.  Let  $\Phi(N,E)=(\omega_1,\ldots,\omega_n)\in\mathcal{M}_{n,1}$ be the sequence of paths defined by
\begin{itemize}
  \item $\omega_{i}=E^{a_{i}}F^{i-1-a_{i}}$, for any $i$ such that $1\leq i \leq n$ and $i\not=k-1, k$.  This path joins $(0,i)$ to $(i-1,i)$.
  \item $\omega_{k-1}=E^{a_{k-1}}F^{\beta}NF^{k-2-a_{k-1}-\beta}$.  This path joins $(0,k-1)$ to $(k-1,k)$.
  \item $\omega_{k}=E^{a_{k}}SE^{b}F^{k-2-a_{k}-b}$. This path joins $(0,k)$ to $(k-2,k-1)$.
\end{itemize}
\end{dfn}
 
{\bf Remark.\,\,} Lemma~\ref{char} shows that $\Phi(N,E)\in\mathcal{M}_{n,1}$.
 

\vspace{2ex}

{\bf Example.\,\,}  The generalized inversion table $(10;\,\, 0,0,2,2,0,0,1,5,0,3,6,6;\,\, 4,5)$ of $(N_{0},3)$ (from figure~\ref{fig2}) is mapped to the second mixed configuration of figure~\ref{fig6}.

\vspace{2ex}

\begin{thm}  The function $\Phi$ is a bijection from
$\mathcal{N}_{n,1}$ to $\mathcal{M}_{n,1}$.
\end{thm}

{\bf Proof.\,\,}  We already know that $\Phi$  maps $\mathcal{N}_{n,1}$ to $\mathcal{M}_{n,1}$. Let $(\omega_1,\ldots,\omega_n) \in \mathcal{M}_{n,1}$.  The paths can be described as sequences of steps as in definition~\ref{defphi} encoding a unique sequence of non-negative integers  $(k;a_{1},\ldots,a_{n};b,\beta)$. Since the paths are in the grid $\mathcal{G}_{n}$, we have $a_{i}\leq i-1$, $a_{k-1}+\beta \leq k-2$ and $a_{k}+b\leq k-2$.  The non-intersecting condition is equivalent to $a_{k-1}<a_{k}$ and $a_{k-1}+\beta<a_{k}+b$.  Clearly  we must have $2\leq k \leq n$.  But if $k=2$ then  $a_{1}<a_{2}\leq  a_{2}+b\leq 0$, a contradiction. Thus $(k;a_{1},\ldots,a_{n};b,\beta)$ is a generalized inversion table of a unique $(N,E)\in \mathcal{N}_{n,1}$. \quad $\Box$

\vspace{2ex}

\begin{figure} 
	\begin{center}
           \includegraphics[scale=1.0]{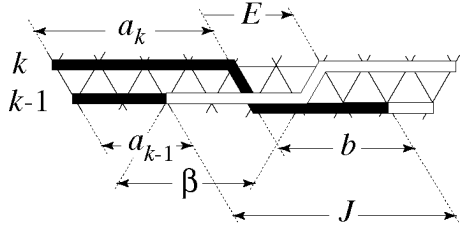} 
	\end{center}

	\caption{ The combinatorial interpretations of the parameters $E$ and  $J$ on mixed configurations. \label{fig7}}
\end{figure}

\begin{thm} \label{interpr}
Let $A\in \mathcal{A}_{n,1}$, $(N,E)= \Lambda(A)$, $(k; a_{1},\ldots,a_{n};b,\beta)$  its generalized inversion table and $\Omega=(\omega_1,\ldots,\omega_n)=\Phi(N,E)$.  Then
\begin{enumerate}
  \item $r(A)=r(N)$ is the number of E-steps of $\Omega$ that are at level $n$ (all occurring in path $\omega_{n}$). \label{uno}
  \item $i(A)=i(N)$ is the total number of E-steps and of N-steps of $\Omega$. \label{duo}
  \item $E(A)=E=a_{k-1}+\beta+1-a_{k}$ is the signed distance from the beginning of the S-step to the end of the N-step of $\Omega$ (see figure~\ref{fig7}). \label{tre}
  \item $B(A)=B(N)-E=b-\beta$.\label{quatro}
  \item $J(A)=J(N)=a_{k}-a_{k-1}+b$ is the (non-signed) distance between the junctions of paths $\omega_{k-1}$ and $\omega_{k}$.\label{cinque}
\end{enumerate}
\end{thm}

{\bf Proof.\,\,}  Statement~\ref{uno} follows from $r(A)=r(N)=a_{n}$.  Statement~\ref{duo}, from the fact that $i(A)=i(N)=a_{1}+\cdots + a_{n}+ b +1$ ($b+1$ is the contribution, after cancellations, of the $-1$ and of the right $1$).  Statement~\ref{tre} is obvious.  As for statement~\ref{quatro}, we have:
$$B(A)=B(N)-E=c(N)-\ell(N)-E=b-(a_{k}-a_{k-1}-1)-E=b-\beta\,.$$
Statement~\ref{cinque} is true, since $J(A)=J(N)=c(N)+\ell(N)+1=b+a_{k}-a_{k-1}$. \quad $\Box$
  
\vspace{2ex}

\begin{cor} \label{exch}
Let $A\in \mathcal{A}_{n,1}$,  $(k; a_{1},\ldots,a_{n};b,\beta)$  the generalized inversion table of $\Lambda(A)$ and $\Omega=(\omega_1,\ldots,\omega_n)=\Phi(\Lambda(A))$.  Then
$(k; a_{1},\ldots,a_{n};b,\beta')$, where $\beta'=b-\beta +a_{k}-a_{k-1} -1$, is the generalized inversion table of $\Lambda(A')$.
\end{cor}

{\bf Proof.\,\,} The generalized inversion table of $\Lambda(A')=(N,E(A'))=(N,B(A))$ is clearly $(k; a_{1},\ldots,a_{n};b,\beta')$ where $\beta'=B(A)+\ell(N)=b-\beta +a_{k}-a_{k-1} -1$.  \quad $\Box$

\vspace{2ex}

Thus the involution $(\,\,)'$ on $\mathcal{A}_{n,1}$ translates into the involution $\Phi \circ \xi \circ \Phi^{-1}$  on $\mathcal{M}_{n,1}$ where it take a very simple form:  $\Phi(\Lambda(A'))$ is obtained from $\Phi(\Lambda(A))$ by replacing $\beta$ by $\beta'$.  The replacement affects  only the Right part of one of the paths.  
 
\vspace{3ex}

\section{Duality and Mixed Configurations}

The counterpart of vertical reflection of alternating sign matrices is an even more useful involution.  As in the case of TB-configuration, this symmetry of mixed configurations will also be a variant of Gessel-Viennot 
paths duality.  

\vspace{2ex}

\begin{figure} 
	\begin{center}
           \includegraphics[scale=1.0]{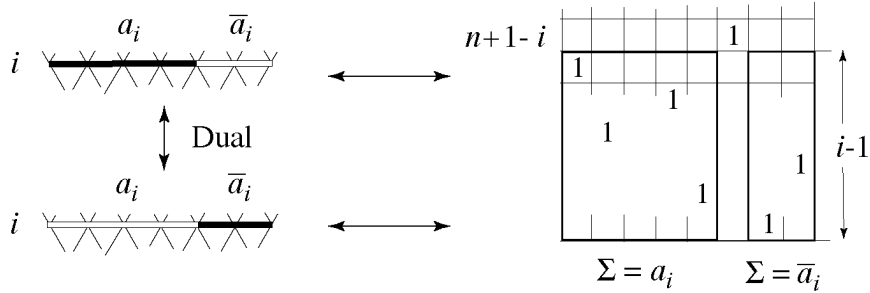} 
	\end{center}

	\caption{Duality on $\mathcal{M}_{n,0}$ (before the final reflection) and its relation with vertical reflection on permutation matrices. \label{fig8}}
\end{figure}

\vspace{2ex}

First, we examine paths duality for mixed configurations $\Omega=(\omega_{1},\ldots,\omega_{n})\in \mathcal{M}_{n,0}$.  We saw how to extract the inversion table  $(a_{i})_{i=1}^{n}$, which corresponds to a unique permutation matrix $P$.  The dual $\overline{\Omega}$ of $\Omega$ is obtained by ``complementing'' the Left and the Right parts (separately) of each path, leading to the sequence $(\overline{a}_{i})_{i=1}^{n}=(i-1-a_{i})_{i=1}^{n}$ which is clearly the inversion table of $\overline{P}$.  In fact, complementation is characterized by the statement:
$$(i,\ell) \text{ is a junction of }\Omega \quad \text{ iff } \quad 
   (\ell-i,\ell) \text{ is a junction of }\overline{\Omega}\,.$$
Therefore, complementation preserves (before the final reflection) the positions of the junctions (see figure~\ref{fig8}).

\vspace{2ex}

For mixed configurations $\Omega=(\omega_{1},\ldots,\omega_{n})\in \mathcal{M}_{n,1}$ (or more generally $\mathcal{M}_{n,s}$), the procedure is similar:
\begin{itemize}
  \item $(i,\ell)$  is a junction of $\Omega$ \quad  iff  \quad 
   $(\ell-i,\ell)$  is a junction of $\overline{\Omega}$,
   
  \item $(i,\ell)$  is the starting vertex of a S step of $\Omega$ \quad iff  \quad $(\ell-i,\ell)$  is the starting vertex of a S step of $\overline{\Omega}$,
  
   \item $(i,\ell)$  is the ending vertex of a N step of $\Omega$ \quad iff  \quad  $(\ell-i,\ell)$ is the ending vertex of a N step of $\overline{\Omega}$.
\end{itemize}
This procedure corresponds to Gessel-Viennot paths duality of the Left and the Right parts of $\Omega$ (separately).  

\begin{figure} 
	\begin{center}
           \includegraphics[scale=1.0]{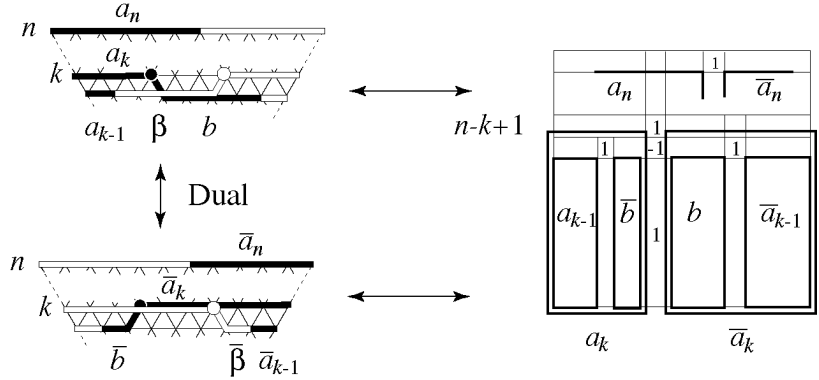} 
	\end{center}

	\caption{Duality on $\mathcal{M}_{n,1}$ (before the final reflection) and its relation with vertical reflection of neutral alternating sign matrices.\label{fig9}}
\end{figure}

\vspace{2ex}
Figure~\ref{fig9} (left) shows the result.   Observe that duality (before the final vertical reflection) preserves the positions of the junctions and of the highest vertices of the S-step and of the  N-steps.  These observations suffice to prove the following lemma (see figure~\ref{fig9} (right)).

\begin{lem}\label{dualite}
Let $\Omega\in \mathcal{M}_{n,1}$, with dual $\overline{\Omega}\in \mathcal{M}_{n,1}$, encoding (respectively) the generalized inversion tables $(k; a_{1},\ldots,a_{n};b,\beta)$ and  $(k; \overline{a}_{1},\ldots,\overline{a}_{n};\overline{b},\overline{\beta})$.   Then
\begin{eqnarray*} 
    \overline{a}_{i} & = & i-1-a_{i} \qquad \mbox{for $i\not=k-1$}, \\
    \overline{a}_{k-1} & = & k-2-a_{k}-b, \\
    \overline{b} & = & a_{k}-1-a_{k-1}, \\
    \overline{\beta} & = & a_{k}+b-a_{k-1}-\beta -1\,. 
\end{eqnarray*}
\end{lem}

\vspace{2ex}

Notice that $\overline{\beta}$ is the same as $\beta'$ of corollary~\ref{exch}.

\begin{thm} \label{dualAM}
Let $A\in \mathcal{A}_{n,1}$ and $\Omega=\Phi(\Lambda(A))$.  Then $\overline{\Omega}= \Phi(\Lambda(\overline{A}))$.
\end{thm}

{\bf Proof.\,\,} It suffices to prove that the generalized inversion tables 
$(k; a_{1},\ldots,a_{n};b,\beta)$ of $\Lambda(A)$ and  $(k;\overline{a}_{1},\ldots,\overline{a}_{n};\overline{b},\overline{\beta})$ of $\Lambda(\overline{A})$ are related as in lemma~\ref{dualite}.  Let 
$(N,E)=\Lambda(A)$ so that $(\overline{N},-E)=\Lambda(\overline{A})$.
 We have:
\begin{itemize}
  \item $a_{i}+\overline{a_{i}}=i-1$ (for $i\not=k-1$) since the sum of the entries of $N$ below row $n+1-i$ is $i-1$.
  \item $a_{k}+b+\overline{a}_{k-1}=k-2$ since the left hand side is the sum of the entries below the closing row (in position $n+2-k$). (The leading 1 is replaced by the lowest 1 in the opening column
  .
  \item $\overline{b}=c(\overline{N})=\ell(N)=a_{k}-a_{k-1}-1$.
  \item  Since  $\beta = E+\ell(N)$, we have:   
  \begin{eqnarray*}
              \overline{\beta} & = & -E+\ell(\overline{N}) \\
                                      & = & -\beta +\ell(N)+c(N)   \\
                                      & = & -\beta +a_{k}-a_{k-1}-1+b\,. \quad \Box
          \end{eqnarray*}
\end{itemize}                 
                                                                                
\end{document}